\documentstyle[epsf,12pt]{article}
\def\Bbb R{{\rm \bf R}}
\def\proclaim#1{\vskip2mm{\bf #1}\em}
\def\endproclaim{\em \vskip2mm}
\def\tag#1{\eqno(#1)}
\def\gathered{\begin{array}{c}}
\def\endgathered{\end{array}}
\def\text{\mbox}

\begin{document}

\title {Travel Time and Heat Equation.\\
One space dimensional case II}
\author{Masaru IKEHATA\\
Department of Mathematics,
Faculty of Engineering\\
Gunma University, Kiryu 376-8515, JAPAN}
\maketitle
\begin{abstract}
Three inverse boundary value problems for the heat equations in
one space dimension are considered. Those three problems are:
extracting an unknown interface in a heat conductive material, an
unknown boundary in a layered material or a material with a smooth
heat conductivity by employing a single set of the temperature and
heat flux on a known boundary as the observation data.  Some
extraction formulae of those discontinuities which suggest a
relationship between the travel time of a {\it virtual} signal and the
observation data are given by applying the enclosure method to the
problems.

\noindent
AMS: 35R30, 80A23

\end{abstract}

\section{Introduction}
In \cite{I3} we considered an
inverse boundary value problem for the heat equation in one space
dimension. First we recall the problem. Let $a>0$ and $T>0$.
Let $u=u(x,t)$
be a solution of the problem:
$$\begin{array}{c}
\displaystyle
u_t=u_{xx}\,\,\text{in}\,]0,\,a[\times]0,\,T[,\\
\\
\displaystyle
u_x(a,t)=0\,\,\text{for}\,t\in\,]0,\,T[,\\
\\
\displaystyle
u(x,0)=0\,\,\text{in}\,]0,\,a[.
\end{array}
\tag {1.1}
$$

\noindent Then the problem considered therein is: extract $a$ from
a {\it single set} of the data $u(0,t)$ and $u_x(0,t)$ for $0<t<T$.

\noindent This is a simplest one space dimension version of the
problem of domain determination which is a typical inverse
boundary value problem for the heat equation and related to the
thermal imaging of unknown {\it discontinuity} such as cavity,
defect or inclusion inside a heat conductive body. There are
extensive studies for the uniqueness and stability issues of this
type of problems in multi dimensions. See \cite{BC, CRV} and
references therein for several results on the issues. However, in
our opinion, seeking an analytical formula that directly connects
information about discontinuity with the data also yields another
view for understanding of the problems.

Recently, some
new analytical methods for the such type inverse problems that the
governing equations are elliptic were introduced. In particular,
the {\it probe method} (\cite{I0, I4}) and {\it factorization
method} (\cite{K, K2}) gave ways of extracting
unknown discontinuity from the data that are given by
the Dirichlet-to-Neumann map in inverse boundary value problems;
the far field operator or the restriction of the
scattered fields onto a sphere surrounding unknown discontinuity
that are exerted by infinitely many point sources located on the
sphere in inverse obstacle scattering problems.

\noindent
These methods require infinitely many data.
By the way, the {\it enclosure method} also gives a way of extracting
the {\it convex hull} of unknown discontinuity from the
Dirichlet-to-Neuman map (\cite{Ie}). However in some cases, using
the idea of the enclosure method, one can give a way of extracting
unknown discontinuity by a single set of the Dirichlet and Neumann
data (see \cite{I1, Ie2}).  In particular, the result in
\cite{Ie2} gave a constructive proof of a uniqueness theorem
established in \cite{FI} and the numerical implementation of a
reconstruction algorithm of the convex hull of unknown polygonal
inclusion or cavity has been done in \cite{IO, IO2}. Thus it is
quite interesting whether the enclosure method still works for the
inverse problem for the heat equation mentioned above. In
\cite{I3} we found a simple extraction formula of $a$ by using the
idea of the enclosure method. Needless to say, this result gives a
new, constructive and simple proof of the uniqueness theorem: the
data $u(0,t)$ and $u_x(0,t)$ for $0<t<T$ uniquely determine $a$
(under a suitable condition on $u_x(0,t)$). However, it should be
pointed out that the result in \cite{I3} is against the past
experiences in the study of inverse boundary value problems for
elliptic equations. Here we review the result and point out the
difference.

\noindent
Let $c$ be an arbitrary positive number.
Set
$$\displaystyle
z=-c\tau\left(1+i\sqrt{1-\frac{1}{c^2\tau}}\,\right),\,\,\tau>c^{-2}.
\tag {1.2}
$$
Let
$$\displaystyle
v(x,t)=e^{-z^2t}\,e^{xz}.
\tag {1.3}
$$
\noindent
Given $s\in\Bbb R$ we introduce the {\it indicator function}
$I_{c}(\tau;s)$:
$$\displaystyle
I_{c}(\tau;s)
=e^{\tau s}
\int_0^T\left(-v_x(0,t)u(0,t)+u_x(0,t)v(0,t)\right)dt,\,\,\tau>c^{-2}
$$
where $u$ satisfies (1.1) and $v$ is the function given by (1.3).

\noindent
Assume that we know a positive
number $M$ such that $M\ge a$. Let $c$ be an arbitrary positive
number satisfying $2Mc<T$. Assume that $u_x(0,t)=1$(note that this
is just for simplicity of description). Then we have the formula
$$\displaystyle
\lim_{\tau\longrightarrow\infty}
\frac{\displaystyle\log\vert I_{c}(\tau;0)\vert}
{\tau}
=-2ca.
\tag {1.4}
$$
and the following statements are true:

if $s\le 2ca$, then
$\displaystyle\lim_{\tau\longrightarrow\infty}\vert I_{c}(\tau;s)\vert=0$;

if $s>2ca$, then
$\displaystyle\lim_{\tau\longrightarrow\infty}\vert I_{c}(\tau;s)\vert=\infty$.

\noindent
Why is this interesting?  The reason consists of two points.

The first point is an unexpected asymptotic behaviour of the
indicator function. More precisely,
integration by parts gives
$$
\displaystyle
I_{c}(\tau;s)
=-e^{\tau s}\int_0^Tu(a,t)v_x(a,t)dt+O(e^{-\tau(T-s)}).
$$
The function $e^{\tau s}v$ has the special character

\noindent
$\bullet$ if $s<cx+t$, then $\lim_{\tau\longrightarrow\infty}e^{\tau s}\vert v(x,t)\vert=0$

\noindent
$\bullet$ if $s>cx+t$, then $\lim_{\tau\longrightarrow\infty}e^{\tau s}\vert v(x,t)\vert=\infty$.

\noindent Therefore if $T>ca>s$, then $\vert I_c(\tau;s)\vert$ is
exponentially decaying as $\tau\longrightarrow\infty$.  Since
$e^{\tau s}\vert v_x(a,t)\vert$ is exponentially growing in the
case when $T>s>ca$ and $0<t<s-ca$, usually we expect $\vert
I_c(\tau;s)\vert$ is exponentially growing as
$\tau\longrightarrow\infty$. Then a past experience suggests that
the right hand side of (1.4) should give $-ca$ which has the
meaning:
$$
-ca=\sup\{\left(\begin{array}{c}\displaystyle x
\\
\displaystyle
t\end{array}\right)\cdot
\left(\begin{array}{c}\displaystyle -c\\
\displaystyle
-1\end{array}\right)\,\vert\,(x,t)\in\,]a,\,\infty[\,\times\,]0,\,T[\}.
$$
This is nothing but the value of the support function for the unknown
domain $]a,\,\infty[\,\times\,]0,\,T[$ at the direction $(-c,\,-1)^T$.
However, in fact, we obtained $-2ca$ which has the meaning:
$$
-2ca=\sup\{\left(\begin{array}{c}\displaystyle x\\
\displaystyle
t\end{array}\right)\cdot
\left(\begin{array}{c} \displaystyle
-c\\
\displaystyle
-1\end{array}
\right)\,\vert\,(x,t)\in\,]2a,\,\infty[\,\times\,]0,\,T[\}.
$$
Since the $u$ of (1.1) can be extended as a solution of the heat
equation onto $]a,\,2a[$ by the reflection $x\longmapsto 2a-x$,
one may think that this is because of the simple Neumann boundary
condition at $x=a$.  However, in \cite{I3} we saw the same
phenomenon to the case of the Robin boundary condition at $x=a$.
So we guess that this is a {\it universal phenomenon}.  Based on this
belief, therein we gave another interpretation of
$2ca$. It is the {\it travel time} of a {\it virtual} signal with an arbitrary
fixed propagation speed $1/c$ that starts at the known boundary
$x=0$ and the initial time $t=0$, reflects at another unknown
boundary $x=a$ and returns to the original boundary $x=0$.

The second point is: needless to say, the formula (1.4) yields a new, simple, constructive proof
of the uniqueness theorem.

$\bullet$  Let
$T$ be a fixed arbitrary positive number.
For $j=1,2$ let $u_j$ satisfy (1.1) with $a=a_j(>0)$ and
$(u_j)_x(0,t)=1$ for all $t\in\,]0,\,T[$.
If $u_1(0,t)=u_2(0,t)$ for all $t\in\,]0,\,T[$, then
$a_1=a_2$.

\noindent
A standard and traditional proof of this type uniqueness theorem
(see, e.g., \cite{BC})
is done by using a contradiction argument
and starts with assuming, say $a_1<a_2$.
Then the uniqueness of the lateral Cauchy problem for the heat equation
gives $u_1(x,t)=u_2(x,t)$ for $(x,t)\in\,]0,\,a_1[\,\times\,]0,\,T[$.  This yields that $u_2$ satisfies
$(u_2)_x(a_1,t)=(u_1)_x(a_1,t)=0$.  Since $u_2$ satisfies the heat equation in $]a_1,\,a_2[\,\times\,]0,\,T[$,
the Neumann boundary conditions at $x=a_1,\,a_2$ and the initial condition $u_2(x,0)=0$, we obtain
$u_2(x,t)=0$ in $]a_1,\,a_2[\,\times\,]0,\,T[$.  Then the unique continuation theorem for the solution
of the heat equation gives $u_2(x,t)=0$ in $]0,\,a_2[\,\times\,]0,\,T[$ and this thus yields
$(u_2)_x(0,t)=0$.  Contradiction.

\noindent Clearly this argument can relax the condition on
$u_x(0,t)$ and tells us the importance of the uniqueness of the
lateral Cauchy problem or the unique continuation theorem for the
heat equation. However, this type proof gives no information about
how to extract unknown $a$ from the data $u(0,t)$ and $u_x(0,t)$.

\noindent
It should be pointed out that, in \cite{BC} another argument
in the case when $T=\infty$ and one
has the data $u(0,t)$ and $u_x(0,t)$ for all $t\in\,]0,\,\infty[$,
is introduced.  Starting with assuming $a_1<a_2$, we see that
$u_2$ satisfies the Neumann boundary conditions at $x=a_1,\,a_2$.
This is same as above.  However, here they do not make use of the
initial condition $u_2(x,0)=0$.  Instead they have the identity
$$\displaystyle
\frac{d}{dt}\int_{a_1}^{a_2}u_2(x,t)dt=\int_{a_1}^{a_2}(u_2)_{xx}(x,t)dt=0.
$$
One the other hand, using the eigenfunction expansion, they show that
$$\displaystyle
\lim_{t\longrightarrow\infty}\vert\int_{a_1}^{a_2}u_2(x,t)dt\vert=\infty.
$$
Contradiction. The proof is quite interesting, however, clearly in
this argument the assumption $T=\infty$ is essential and the same
comment works also for this proof.

\noindent Summing up, one should find another type proof that
tells us the information about how to extract unknown $a$ from the
data $u(0,t)$ and $u_x(0,t)$ for $t\in\,]0,\,T[$ with $T<\infty$.
We think that our proof presented in this paper gives an answer to
this natural question.

The aim of this paper is: to confirm further that the interpretation
in the first point still works, at least, in one space dimensional case
by considering
three typical inverse boundary value problems for the heat
equations.  Those three problems are:
extracting an unknown interface in a conductive material,
an unknown boundary in a layered material or a material with a smooth
conductivity.

\noindent In a future study we will consider the multidimensional
version of those problems.

\noindent
{\bf\noindent Remark 1.1.}
In this paper we always consider
the solutions of the heat equations in the context of a
variational formulation.  In particular, every solutions in this paper
belongs to the space $W(0,\,T;H^1(\Omega), (H^1(\Omega))')$ with $u_x(0,t), u_x(a,t)\in L^2(0,T)$
where $\Omega=]0,\,a[$ and satisfies the governing equation in a weak sense.
We refer the reader to \cite{DL} for the detail.

\section{Statement of the results}

{\bf\noindent 2.1. Extracting interface}

Let $0<b<a$.
Define
$$
\gamma(x)=\left\{
\begin{array}{lr}
\displaystyle \gamma_1,&\quad \text{if $0<x<b$,}\\
\\
\displaystyle \gamma_2,&\quad \text{if $b<x<a$}
\end{array}
\right.
$$
where both $\gamma_1$ and $\gamma_2$ are positive constants and satisfies $\gamma_2\not=\gamma_1$.

Let $u$ be an arbitrary solution of the problem:
$$\begin{array}{c}
\displaystyle
u_t=(\gamma\,u_{x})_{x}\,\,\text{in}\,]0,\,a[\times\,]0,\,T[,\\
\\
\displaystyle
u(x,0)=0\,\text{in}\,\,]0,\,a[.
\end{array}
\tag {2.1}
$$
We assume that both $\gamma_1$ is {\it known}
and that $a$, $b$ and $\gamma_2$ are all {\it unknown}.

{\bf\noindent Inverse Problem A.} Extract $b$ from $u(0,t)$ and
$\gamma_1\,u_x(0,t)$ for $0<t<T$.

\noindent
Let $c$ be an arbitrary positive number. Let
$$\displaystyle
v(x,t)=e^{-z^2t}\Psi(x,z)
\tag {2.2}
$$
where $\Psi(x,z)=e^{x\,z_1}$, $z_1=z/\sqrt{\gamma_1}$ and $z$ is given by (1.2).

\noindent
The function $v$ is a complex valued function and satisfies the backward heat equation $v_t+\gamma_1\,v_{xx}=0$
in the whole space-time.

{\bf\noindent Definition 2.1.}
Given $c>0$
define the {\it indicator function}
$I_{c}(\tau)$ by the formula
$$\displaystyle
I_{c}(\tau)
=\int_0^T\left(-\gamma_1\,v_x(0,t)u(0,t)+\gamma_1\,u_x(0,t)v(0,t)\right)dt,\,\,\tau>c^{-2}
$$
where $u$ satisfies (2.1) and $v$ is the function given by (2.2).

Define
$$\displaystyle
w(x)=w(x,\tau)=\int_0^T\,e^{-z^2\,t}\,u(x,t)dt,\,\,0<x<a.
$$
Then this $w$ satisfies
$$\displaystyle
(\gamma\,w')'-z^2\,w=e^{-z^2\,T}\,u(x,T)\,\,\text{in}\,]0,\,a[
\tag {2.3}
$$
and
$$
\displaystyle
w'(0)=\int_0^Te^{-z^2\,t}u_x(0,t)dt,\,\,
w'(a)=\int_0^Te^{-z^2\,t}u_x(a,t)dt.
$$

Our first result is the following theorem:

\proclaim{\noindent Theorem 2.1.}
Assume that we know a positive number $M$ such that
$M\ge 2b/\sqrt{\gamma_1}$.  Let $c$ be an arbitrary positive number satisfying
$Mc<T$.
Assume that
$$\displaystyle
\lim_{\tau\longrightarrow\infty}\,\vert\frac{w'(a)}{w'(0)}\vert
\exp\,\left(\displaystyle c\,\tau\,(\frac{b}{\sqrt{\gamma_1}}-\frac{(a-b)}{\sqrt{\gamma_2}})\right)=0
\tag {2.4}
$$
and there exist a positive constant $C$, positive
number $\tau_0$, real number $\mu$ such that, for all
$\tau>\tau_0$
$$
\displaystyle
C\,\tau^{\mu}\,\le\vert w'(0)\vert.
\tag {2.5}
$$
Then the formula
$$\displaystyle
\lim_{\tau\longrightarrow\infty}
\frac{\displaystyle\log\vert I_{c}(\tau)\vert}
{\tau}
=-2\,\frac{c\,b}{\sqrt{\gamma_1}},
$$
is valid.

\endproclaim

\noindent Note that the boundary condition at $x=a$ is not
specified.  However, the condition (2.4) implicitly restricts a
possible boundary condition at $x=a$. The situation dramatically
changes in the case when (a)
$b/\sqrt{\gamma_1}<(a-b)/\sqrt{\gamma_2}$ or (b)
$b/\sqrt{\gamma_1}>(a-b)/\sqrt{\gamma_2}$.

\noindent In case (a) the condition (2.5) automatically ensures
that (2.4) is valid since we have always $w'(a)=O(1)$ as
$\tau\longrightarrow\infty$(Remark 1.1). This is reasonable under
our interpretation since the signal started from $x=a$ at $t=0$
arrives at $x=0$ after the arrival of the signal started at $x=0$
at $t=0$ (see Figure 1).

\begin{figure}
\begin{center}
\epsfxsize=9cm
\epsfysize=9cm
\epsfbox{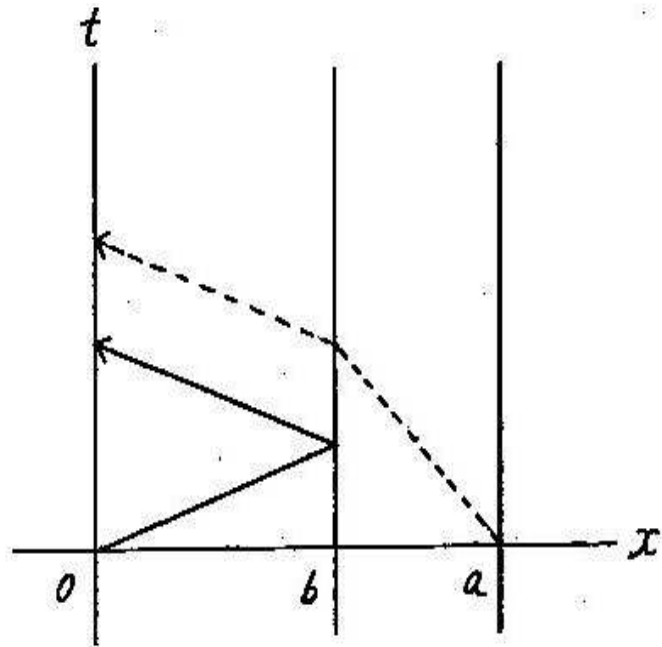}
\caption{(a) $b/\sqrt{\gamma_1}<(a-b)/\sqrt{\gamma_2}$.}
\end{center}
\end{figure}

\noindent However in case (b) to ensure (2.4) $w'(a)$ has to decay
exponentially as $\tau\longrightarrow\infty$.  This is a strong
restriction and can be interpreted as a condition that kills a
signal started from $x=a$ at $t=0$ (see Figure 2).

\begin{figure}
\begin{center}
\epsfxsize=9cm
\epsfysize=9cm
\epsfbox{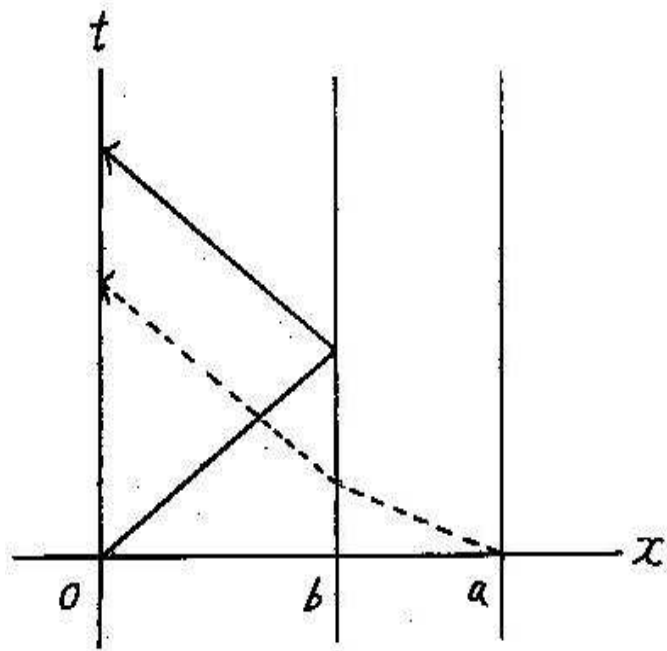}
\caption{(b) $b/\sqrt{\gamma_1}>(a-b)/\sqrt{\gamma_2}$.}
\end{center}
\end{figure}

The condition (2.5) gives a restriction on the behaviour of the flux $u_x(0,t)$ as $t\downarrow 0$.
Let $\delta$ satisfy $0<\delta<T$ and $m\ge 0$ be an integer.
Since
$$\displaystyle
z^2=\tau+i2c^2\tau^2\sqrt{1-\frac{1}{c^2\tau}},
$$
a change of variable yields, as $\tau\longrightarrow\infty$
$$\begin{array}{c}
\displaystyle

\int_0^{\delta}e^{-z^2 t}t^m dt
=\frac{1}{\tau^{m+1}}\int_0^{\tau\delta}e^{\displaystyle-\xi(1+i2c^2\tau\sqrt{1-\frac{1}{c^2\tau}}\,)}\xi^md\xi\\
\\
\displaystyle
=\frac{1}{\tau^{m+1}}\int_0^{\infty}e^{\displaystyle-c^2\tau(1+i\sqrt{1-\frac{1}{c^2\tau}}\,)^2\xi}\xi^md\xi
+O(\tau^{-\infty})
\\
\\
\displaystyle
=\frac{1}{\tau^{2(m+1)}}
\int_0^{\infty}e^{\displaystyle-c^2(1+i\sqrt{1-\frac{1}{c^2\tau}}\,)^2 t}t^mdt
+O(\tau^{-\infty}).
\end{array}
$$
Integration by parts gives
$$
\displaystyle
\int_0^{\infty}e^{\displaystyle-c^2(1+i\sqrt{1-\frac{1}{c^2\tau}}\,)^2
t}t^mdt =m!K(\tau)^{m+1}
$$
where
$$\displaystyle
K(\tau)=\frac{1}{\displaystyle
c^2(1+i\sqrt{1-\frac{1}{c^2\tau}})^2} \longrightarrow
\frac{-i}{2c^2}.
$$
Therefore we obtain, as $\tau\longrightarrow\infty$
$$
\displaystyle
\tau^{2(m+1)}\int_0^{\delta}e^{-z^2 t}t^mdt
\longrightarrow m!(-\frac{i}{2c^2})^{m+1}\not=0.
$$
It is easy to see that if $f\in L^2(0,\,T)$ and for a suitable positive constant $C$,
$\vert f(t)\vert\le Ct^{m'}$ a. e. in $]0,\,\delta[$, then,
as $\tau\longrightarrow\infty$
$$\displaystyle
\int_0^Te^{-z^2 t}f(t)dt=O(\tau^{-(m'+1)}).
$$
Now assume that, we have, for some $\delta$, a positive constant $C$
$$\displaystyle
\vert u_x(0,t)-(g_lt^l+g_{l+1}t^{l+1}+\cdots+g_nt^n)\vert\le Ct^{m'+1},\,\,\text{a. e. in $]0,\,\delta[$}
$$
where $n$, $l$ and $m'$ are integers and and satisfy $n\ge l\ge 0$ and $m'>2l+1$; $g_l$,
... $g_n$ are constants and satisfy $g_l\not=0$.  Then from the
computation above we see that, as $\tau\longrightarrow\infty$
$$\displaystyle
\tau^{2(l+1)}\int_0^{T}e^{-z^2 t}u_x(0,t)dt
\longrightarrow g_ll!(-\frac{i}{2c^2})^{l+1}.
$$
This implies that $u_x(0,t)$ satisfies (2.5).
In particular, the condition (2.5) is satisfied if
$u_x(0,t)$ satisfies one of the following conditions

$\bullet$  $u_x(0,t)\in C^{2}([0,\,\delta])$
and $u_x(0,0)\not=0$

$\bullet$  $u_x(0,t)\in C^{2l+2}([0,\,\delta])$ with $l\ge 1$
and $t=0$ is the zero point of $u_x(0,t)$ with order $l$.

{\bf\noindent 2.2.  Extracting unknown boundary. Layered material}

Let $0=b_0<b_1<b_2<\cdots<b_m=a$.
$$
\gamma(x)=\left\{
\begin{array}{lr}
\displaystyle\gamma_1,&\quad \text{if $b_0<x<b_1$,}\\
\\
\displaystyle\gamma_2,&\quad \text{if $b_1<x<b_2$,}\\
\\
\vdots\\
\\
\displaystyle\gamma_m, &\quad\text{if $b_{m-1}<x<b_m$}
\end{array}
\right.
$$
where $\displaystyle\gamma_1,\gamma_2,\cdots,\gamma_m$ are positive constants and $m\ge 1$.

\noindent  Let $a>0$. Let $u=u(x,t)$ be an arbitrary solution of
the problem:
$$\begin{array}{c}
\displaystyle
u_t=(\gamma\, u_{x})_{x}\,\,\text{in}\,]0,\,a[\times]0,\,T[,\\
\\
\displaystyle
\gamma_m\,u_x(a,t)+\rho\,u(a,t)=0\,\,\text{for}\,t\in\,]0,\,T[,\\
\\
\displaystyle
u(x,0)=0\,\,\text{in}\,]0,\,a[
\end{array}
\tag {2.6}
$$
here $\rho\ge 0$ is an arbitrary fixed constant.

{\bf\noindent Inverse Problem B.}  Assume that $\gamma_1$, $\cdots$,
$\gamma_m$ are {\it known} and that both of $\rho$ and $a$ are
{\it unknown}.  Extract $a$ from $u(0,t)$ and $\gamma_1\,u_x(0,t)$
for $0<t<T$.

In this subsection, instead of $\Psi$ in the previous subsection the function $\Psi$ in
the following which is a consequence of Lemma 5.1 in Section 5 plays the same role;
$z$ is given by (1.2) and set $z_j=z/\sqrt{\gamma_j}$, $j=1,\cdots,m$.

\proclaim{\noindent Proposition 2.2.}
There exists a positive number
$C=C(c,\gamma_1,\cdots,\gamma_m)\,(>c^{-2})$ such that: for each $\tau>C$
there exist a unique $U=(B_1, A_2, B_2,\cdots, A_{m-1},B_{m-1}, A_m)$
such that the function $\Psi$ defined by the formula
$$
\Psi(x)=\Psi(x;c,\gamma_1,\cdots,\gamma_m,\tau)=\left\{
\begin{array}{lr}
\displaystyle e^{xz_1}+B_1e^{-xz_1},&\quad \text{if $x<b_1$,}\\
\\
\displaystyle A_2e^{xz_2}+B_2e^{-xz_2},&\quad \text{if $b_1<x<b_2$,}\\
\\
\vdots\\
\\
\displaystyle A_{m-1}e^{xz_{m-1}}+B_{m-1}e^{-xz_{m-1}}, &\quad\text{if $b_{m-2}<x<b_{m-1}$,}\\
\\
\displaystyle A_m e^{xz_m}, &\quad\text{if $b_{m-1}<x$}
\end{array}
\right.
$$
satisfies the equation
$$
\displaystyle (\tilde{\gamma}\Psi')'-z^2\Psi=0\,\,\text{in}\,\Bbb R
$$
where
$$
\tilde{\gamma}(x)=\left\{
\begin{array}{lr}
\displaystyle\gamma_1,&\quad \text{if $-\infty<x<b_1$,}\\
\\
\displaystyle\gamma_2,&\quad \text{if $b_1<x<b_2$,}\\
\\
\vdots\\
\\
\displaystyle\gamma_m, &\quad\text{if $b_{m-1}<x<\infty$.}
\end{array}
\right.
$$

\endproclaim

Using $\Psi$ in Proposition 2.2, we define the special solution of the
backward heat equation $v_t+(\gamma v_x)_x=0$ in $\Bbb R$ by the
formula
$$
\displaystyle
v(x,t)=e^{-z^2t}\Psi(x;c,\gamma_1,\cdots,\gamma_m,\tau),\,\,
\tau>C(c;\gamma_1,\cdots,\gamma_m)
$$

\noindent Now we can define the indicator function.

{\bf\noindent Definition 2.2.} Given $c>0$ define the {\it
indicator function} $I_{c}(\tau)$ by the formula
$$\displaystyle
I_{c}(\tau)
=\int_0^T\left(-\gamma_1\,v_x(0,t)u(0,t)+\gamma_1\,u_x(0,t)v(0,t)\right)dt,\,\,
\tau>C(c;\gamma_1,\cdots,\gamma_m)
$$
where $u$ satisfies (2.6).

The following gives a solution to the problem mentioned above and
generalizes Theorem 4.1 of \cite{I3}.

\proclaim{\noindent Theorem 2.3.} Assume that we know a positive
number $M$ such that
$$\displaystyle
M\ge 2
\left(\frac{b_1}{\sqrt{\gamma_1}}+\sum_{j=1}^{m-1}\frac{(b_{j+1}-b_j)}{\sqrt{\gamma_{j+1}}}
\right).
$$
Let $c$ be an arbitrary positive
number satisfying $Mc<T$.  Assume that there exist positive a constant $C$, positive
number $\tau_0(>c^{-2})$, real number $\mu$ such that, for all
$\tau>\tau_0$
$$
\displaystyle
C\,\tau^{\mu}\,\le\vert\int_0^T e^{-z^2\,t}u_x(0,t)dt\vert.
\tag {2.7}
$$

\noindent
Then the formula
$$\displaystyle
\lim_{\tau\longrightarrow\infty}
\frac{\displaystyle\log\vert I_{c}(\tau)\vert}
{\tau}
=-2c\left(\frac{b_1}{\sqrt{\gamma_1}}+\sum_{j=1}^{m-1}\frac{(b_{j+1}-b_j)}{\sqrt{\gamma_{j+1}}}\right),
\tag {2.8}
$$
is valid.

\endproclaim

\noindent
Needless to say, the quantity
$$\displaystyle
2c\left(\frac{b_1}{\sqrt{\gamma_1}}+\sum_{j=1}^{m-1}\frac{(b_{j+1}-b_j)}{\sqrt{\gamma_{j+1}}}\right)
$$
can be interpreted as the travel time of a virtual signal
with propagation speeds
\newline{$\sqrt{\gamma_1}/c, \sqrt{\gamma_2}/c,\cdots,\sqrt{\gamma_m}/c$}
in the layeres $0<x<b_1, b_1<x<b_2,\cdots,b_{m-1}<x<a$, respectively that
starts at the known boundary
$x=0$ and the initial time $t=0$, reflects at another unknown
boundary $x=a$ and returns to the original boundary $x=0$.

{\bf\noindent 2.3.  Extracting unknown boundary. Material with
smooth conductivity}

Let $a>0$.  Let $M\ge a$.
Let $\gamma\in C^2([0,\,M])$ and satisfy $\gamma(x)>0$ for all $x\in\,[0,\,M]$.
\noindent  Let $u=u(x,t)$ be an arbitrary solution of
the problem:
$$\begin{array}{c}
\displaystyle
u_t=(\gamma\, u_{x})_{x}\,\,\text{in}\,]0,\,a[\times]0,\,T[,\\
\\
\displaystyle
\gamma(a)\,u_x(a,t)+\rho\,u(a,t)=0\,\,\text{for}\,t\in\,]0,\,T[,\\
\\
\displaystyle
u(x,0)=0\,\,\text{in}\,]0,\,a[
\end{array}
\tag {2.9}
$$
here $\rho\ge 0$ is an arbitrary fixed constant.

{\bf\noindent Inverse Problem C.}
Assume that both $M$ and $\gamma$ are {\it known} and that both of $\rho$ and $a$ are
{\it unknown}.  Extract $a$ from $u(0,t)$ and $\gamma(0)\,u_x(0,t)$
for $0<t<T$.

We start with a fact which can be deduced from a combination of Theorem 1 of p. 48 in \cite{N}
and the Liouville transform.

\proclaim{\noindent Proposition 2.4.}
Given $z$ with $\text{Re}\,z\le 0$ there exists a solution
$\Psi=\Psi(\,\cdot\,;z,M)$ of the equation $(\gamma y')'-z^2 y=0, \,0<x<M$ such that,
as $\vert z\vert\longrightarrow\infty$
$$\displaystyle
\Psi(x;z,M)=\{\gamma(x)\}^{-1/4}\,
\exp\,\left(z\,\int_0^x\frac{dx}{\sqrt{\gamma(x)}}\right)
\,\{1+O\left(\frac{1}{\vert z\vert}\right)\},
$$
$$\displaystyle
\Psi'(x;z,M)
=z\,\{\gamma(x)\}^{-3/4}\,
\exp\,\left(z\,\int_0^x\frac{dx}{\sqrt{\gamma(x)}}\right)
\,\{1+O\left(\frac{1}{\vert z\vert}\right)\}
$$
uniformly in $x\in\,[0,\,M]$.

\endproclaim

Using $\Psi$ in Proposition 2.4, we define the special solution of the
backward heat equation $v_t+(\gamma v_x)_x=0$ in $]0,\,M[\times\,]0,\,T[$ by the
formula
$$
\displaystyle
v(x,t)=e^{-z^2t}\Psi(x;z,M).
$$

\noindent Now we can define the indicator function.

{\bf\noindent Definition 2.3.} Define the {\it
indicator function} $I(z)$ by the formula
$$\displaystyle
I(z)
=\int_0^T\left(-\gamma(0)\,v_x(0,t)u(0,t)+\gamma(0)\,u_x(0,t)v(0,t)\right)dt,\,\,
\text{Re}\,z\le 0
$$
where $u$ satisfies (2.9).

The following gives two solutions to the problem mentioned above and
generalizes Theorem 4.1 of \cite{I3}.

\proclaim{\noindent Theorem 2.5.} Assume that we know a positive
number $M$ such that $M\ge a$.

(1)
Let $c$ be an arbitrary positive
number satisfying
$$\displaystyle
T>2c\int_0^M\frac{dx}{\sqrt{\gamma(x)}}.
\tag {2.10}
$$
Let $z$ be the number given by (1.2).
Assume that there exist positive constant $C$ and a positive
number $\tau_0(>c^{-2})$, real number $\mu$ such that, for all
$\tau>\tau_0$
$$
\displaystyle
C\,\tau^{\mu}\,\le\vert\int_0^T\,e^{-z^2\, t}\,u_x(0,t)dt\vert.
\tag {2.11}
$$

\noindent
Then the formula
$$\displaystyle
\lim_{\tau\longrightarrow\infty}
\frac{\displaystyle\log\vert I(z)\vert}{\tau}
=-2c\int_0^a\frac{dx}{\sqrt{\gamma(x)}},
\tag {2.12}
$$
is valid.

(2)
Assume that there exist positive constant $C$ and a positive
number $\tau_0$, real number $\mu$ such that, for all
$\tau>\tau_0$
$$
\displaystyle
C\,\tau^{\mu}\,\le\vert\int_0^T\,e^{-\tau^2\, t}\,u_x(0,t)dt\vert.
\tag {2.13}
$$
Then the formula
$$\displaystyle
\lim_{\tau\longrightarrow\infty}
\frac{\displaystyle\log\vert I(-\tau)\vert}{\tau}
=-2\int_0^a\frac{dx}{\sqrt{\gamma(x)}},
\tag {2.14}
$$
is valid.

\endproclaim

\noindent
The quantity
$$\displaystyle
2c\int_0^a\frac{dx}{\sqrt{\gamma(x)}}
$$
coincides with the travel time of a virtual signal with variable propagation speed
$\displaystyle\sqrt{\gamma(x)}/c$ that
starts at the known boundary
$x=0$ and the initial time $t=0$, reflects at another unknown
boundary $x=a$ and returns to the original boundary $x=0$.

The condition (2.13) is less restrictive compared with condition (2.11).
We can easily see that the condition (2.13) is satisfied
if $u_x(0,t)$ satisfies one of the following conditions

$\bullet$  for a positive constant $C$ $u_x(0,t)\ge C$ a.e. in $]0,\,\delta[$

$\bullet$  $u_x(0,t)\in C^{1}([0,\,\delta])$
and $u_x(0,0)\not=0$

$\bullet$  $u_x(0,t)\in C^{l+1}([0,\,\delta])$ with $l\ge 1$
and $t=0$ is the zero point of $u_x(0,t)$ with order $l$.

\noindent
See the end of subsection 2.1 for the comparison.

As a corollary of Theorem 2.5 we obtain a direct proof of {\it uniqueness theorem}:
the data $u(0,t)$ and $u_x(0,t)$ for $0<t<T$ uniquenely determine $a$
provided one of (2.11) for a $c$ satisfying (2.10) or (2.13) is satisfied
and both $M$ and $\gamma$ are known.
Note that $T$ is an arbitrary fixed positive number.
In the proof we never make use of
the uniqueness of the lateral Cauchy problem nor the unique continuation theorem for
the heat equation with a variable coefficient.

\section{Proof of Theorem 2.1}

Since $\Psi$ satisfies $\gamma_1\Psi''-z^2\Psi=0$ in $\Bbb R$,
integration by parts gives the expression of the indicator function:
$$\begin{array}{l}
\displaystyle
I_c(\tau)
=-\gamma_1\,\Psi'(a)w(a)+\gamma_2\,w'(a)\Psi(a)
+(\gamma_1-\gamma_2)\int_b^a w'(x)\Psi'(x)dx\\
\\
\displaystyle
-e^{-z^2\,T}\int_0^au(\xi,T)\Psi(\xi)d\xi.
\end{array}
\tag {3.1}
$$

\noindent
Let $y$ be the solution of the boundary value problem
$$\begin{array}{c}
\displaystyle
(\gamma\,y')'-z^2\,y=0\,\,\text{in}\,]0,\,a[,\\
\\
\displaystyle
y'(0)=w'(0),\,\,y'(a)=w'(a).
\end{array}
$$

\noindent Define the {\it principle part} of the indicator
function
$$\displaystyle
I_c^{0}(\tau)=
-\gamma_1\Psi'(a)y(a)+\gamma_2 y'(a)\Psi(a)
+(\gamma_1-\gamma_2)\,\int_b^a y'(x)\Psi'(x)dx.
\tag {3.2}
$$

It is easy to see that Theorem 2.1 is a direct consequence of the following two lemmas

\proclaim{\noindent Lemma 3.1.}
Assume that $w'(a)$ and $w'(0)$ satisfy (2.4) and (2.5).
Then, choosing suitable positive constants $C'_1$ and $C_2'$ and positive number $\tau_0'$, we have
$$
C_1'\,\tau^{\mu}\le\vert I_c^0(\tau)\vert e^{2\,c\,\tau\,b/\sqrt{\gamma_1}}
\le C_2',\,\,\forall\tau>\tau_0'.
$$

\endproclaim

\proclaim{\noindent Lemma 3.2.}
As $\tau\longrightarrow\infty$
$$\displaystyle
I_c(\tau)=I_c^0(\tau)
+O\left(e^{-\tau\,T}\right).
$$
\endproclaim

First we prove Lemma 3.1.
\noindent
The $y$ has the expression
$$
y(x)=y(x,\tau)=\left\{
\begin{array}{lr}
\displaystyle A_1\,e^{x\,z_1}+B_1\,e^{-x\,z_1},&\quad \text{if $0<x<b$,}\\
\\
\displaystyle A_2\,e^{x\,z_2}+B_2\,e^{-x\,z_2},&\quad \text{if $b<x<a$}
\end{array}
\right.
$$
and has to satisfy the transmission conditions:
$$\displaystyle
y(b-0)=y(b+0),\,\,\gamma_1\,y'(b-0)=\gamma_2\,y'(b+0).
$$
These and the relation $z_2/z_1=\sqrt{\gamma_1/\gamma_2}$ yield the system of equations:
$$\displaystyle
A_1\,e^{b\,z_1}+B_1\,e^{-b\,z_1}
=A_2\,e^{b\,z_2}+B_2\,e^{-b\,z_2};
\tag {3.3}
$$
$$\displaystyle
A_1\,e^{b\,z_1}-B_1\,e^{-b\,z_1}
=\sqrt{\frac{\gamma_2}{\gamma_1}}\,(A_2\,e^{b\,z_2}-B_2\,e^{-b\,z_2}).
\tag {3.4}
$$

\noindent
Moreover the boundary conditions $y'(0)=w'(0)$ and $y'(a)=w'(a)$ yield
$$\displaystyle
z_1\,(A_1-B_1)=w'(0)
\tag {3.5}
$$
and
$$
z_2\,(A_2\,e^{a\,z_2}-B_2\,e^{-a\,z_2})=w'(a).
\tag {3.6}
$$

\noindent
A combination of (3.3) and (3.4) gives
$$\begin{array}{l}
\displaystyle
A_1\,e^{b\,z_1}
=\frac{1}{2}\left(1+\sqrt{\frac{\gamma_2}{\gamma_1}}\,\right)\,A_2\,e^{b\,z_2}
+\frac{1}{2}
\left(1-\sqrt{\frac{\gamma_2}{\gamma_1}}\,\right)\,B_2\,e^{-b\,z_2},\\
\\
\displaystyle
B_1\,e^{-b\,z_1}
=\frac{1}{2}\left(1-\sqrt{\frac{\gamma_2}{\gamma_1}}\,\right)\,A_2\,e^{b\,z_2}
+\frac{1}{2}
\left(1+\sqrt{\frac{\gamma_2}{\gamma_1}}\,\right)
\,B_2\,e^{-b\,z_2}.
\end{array}
\tag {3.7}
$$

\noindent
Define
$$\displaystyle
T_{kl}=\frac{2\,\sqrt{\gamma_k}}{\sqrt{\gamma_k}+\sqrt{\gamma_l}},\,\,
R_{kl}=\frac{\sqrt{\gamma_k}-\sqrt{\gamma_l}}
{\sqrt{\gamma_k}+\sqrt{\gamma_l}}.
$$
We have
$$
\displaystyle
\frac{1}{2}
\left(1+\sqrt{\frac{\gamma_2}{\gamma_1}}\,\right)=
\frac{1}{T_{12}},\,\,
\frac{1}{2}\left(1-\sqrt{\frac{\gamma_2}{\gamma_1}}\,\right)
=\frac{R_{12}}{T_{12}}.
$$
Then (3.7) becomes equations
$$\displaystyle
A_1\,e^{b\,z_1}=\frac{1}{T_{12}}\,A_2\,e^{b\,z_2}
+\frac{R_{12}}{T_{12}}\,B_2\,e^{-b\,z_2},
\tag {3.8}
$$
$$\displaystyle
B_1\,e^{-b\,z_1}=\frac{R_{12}}{T_{12}}
\,A_2\,e^{b\,z_2}
+\frac{1}{T_{12}}\,B_2\,e^{-b\,z_2}.
\tag {3.9}
$$

\noindent
Substituting (3.8) and (3.9) into (3.5) and (3.6), we obtain
$$\displaystyle
R\left(\begin{array}{c} \displaystyle A_2\,e^{b\,z_2}\\
\\
\displaystyle B_2\,e^{-b\,z_2}\end{array}\right)
=\frac{1}{z}\,\left(\begin{array}{c}\displaystyle \sqrt{\gamma_2}\,w'(a)\\
\\
\displaystyle \sqrt{\gamma_1}\,w'(0)\,e^{b\,z_1}\end{array}\right)
\tag {3.10}
$$
where
$$\begin{array}{c}
\displaystyle
R=\left(\begin{array}{c}
\begin{array}{lr} e^{(a-b)\,z_2} & -e^{-(a-b)\,z_2}\end{array}
\\
\\
\displaystyle
\left(\begin{array}{lr} 1 & -e^{2b\,z_1}\end{array}\right)
\left(\begin{array}{lr} \displaystyle \frac{1}{T_{12}} & \displaystyle\frac{R_{12}}{T_{12}}\\
\\
\displaystyle
\frac{R_{12}}{T_{12}} & \displaystyle\frac{1}{T_{12}}\end{array}\right)
\end{array}
\right)\\
\\
\displaystyle
=\left(\begin{array}{lr} \displaystyle e^{(a-b)\,z_2} & \displaystyle -e^{-(a-b)\,z_2}\\
\\
\displaystyle
\frac{1}{T_{12}}\,(1-R_{12}\,e^{2b\,z_1}) & \displaystyle \frac{1}{T_{12}}\,
(R_{12}-e^{2b\,z_1})\end{array}
\right).
\end{array}
$$

\noindent
A direct computation gives
$$\displaystyle
(\text{det}\,R)\,T_{12}\,e^{(a-b)\,z_2}
=1+R_{12}(e^{2(a-b)\,z_2}-e^{2b\,z_1})
-e^{2b\,z_1+2(a-b)\,z_2}.
$$
This yields,
as $\tau\longrightarrow\infty$
$$\displaystyle
(\text{det}\,R)\,T_{12}\,e^{(a-b)\,z_2}
=1+O\left(e^{\displaystyle -2c\tau\,\min\,(b/\sqrt{\gamma_1},\,(a-b)/\sqrt{\gamma_2}\,)}\right).
\tag {3.11}
$$
Therefore $R$ is invertible for sufficiently large $\tau$.
It follows from (3.10) that $A_2$, $B_2$ have the form
$$\displaystyle
A_2\,=\frac{e^{-b\,z_2}}{z(\text{det}\,R)}
\left(\frac{1}{T_{12}}(R_{12}-e^{2b\,z_1})\,\sqrt{\gamma_2}\,w'(a)
+e^{-(a-b)\,z_2}\,\sqrt{\gamma_1}\,w'(0)\,e^{b\,z_1}\right),
\tag {3.12}
$$
$$\displaystyle
B_2\,=
\frac{e^{b\,z_2}}{z(\text{det}\,R)}
\left(\frac{1}{T_{12}}(R_{12}\,e^{2b\,z_1}-1)\,\sqrt{\gamma_2}\,w'(a)
+e^{(a-b)\,z_2}\,\sqrt{\gamma_1}\,w'(0)\,e^{b\,z_1}\right).
\tag {3.13}
$$

\noindent
Using (3.12) and (3.13), we get two crucial formulae:
$$
\displaystyle
(\gamma_1-\gamma_2)\,\int_b^a y'(x)\Psi'(x)dx
=\frac{1}{(\text{det}\,R)\,T_{12}}
\left(C\,\sqrt{\gamma_2}\,w'(a)+D\,T_{12}\,\sqrt{\gamma_1}\,w'(0)\right)
$$
where
$$\begin{array}{l}
\displaystyle
C=(\sqrt{\gamma_1}-\sqrt{\gamma_2})(R_{12}-e^{2b\,z_1})
\,(e^{a\,z_1+(a-b)\,z_2}-e^{b\,z_1})\\
\\
\displaystyle
+(\sqrt{\gamma_1}+\sqrt{\gamma_2})
(R_{12}\,e^{2b\,z_1}-1)\,
(e^{a\,z_1-(a-b)\,z_2}-e^{b\,z_1}),
\\
\\
\displaystyle
D=2\sqrt{\gamma_1}\,e^{(a+b)\,z_1}\\
\\
\displaystyle
-(\sqrt{\gamma_1}\,-\sqrt{\gamma_2})\,e^{2b\,z_1-(a-b)\,z_2}
-(\sqrt{\gamma_1}+\sqrt{\gamma_2})\,
e^{2b\,z_1+(a-b)\,z_2};
\end{array}
$$
$$\displaystyle
-\gamma_1\Psi'(a)y(a)+\gamma_2 y'(a)\Psi(a)
=\frac{1}{(\text{det\,R})\,T_{12}}
\left(\tilde{C}\sqrt{\gamma_2}\,w'(a)+\tilde{D}\,T_{12}\,\sqrt{\gamma_1}\,w'(0)\right)
$$
where
$$\begin{array}{l}
\displaystyle
\tilde{C}=-(\sqrt{\gamma_1}-\sqrt{\gamma_2})(R_{12}-e^{2b\,z_1})
e^{a\,z_1+(a-b)\,z_2}\\
\\
\displaystyle
-(\sqrt{\gamma_1}+\sqrt{\gamma_2})(R_{12}\,e^{2b\,z_1}-1)e^{a\,z_1-(a-b)\,z_2},\\
\\
\displaystyle
\tilde{D}=-2\sqrt{\gamma_1}\,e^{(a+b)\,z_1}.
\end{array}
$$

\noindent
A combination of those and (3.2) gives the representation formula
$$\begin{array}{l}
\displaystyle
(\text{det}\,R)\,T_{12}\,I_c^0(\tau)
=-2\,\sqrt{\gamma_1}\,\sqrt{\gamma_2}\,w'(a)\,e^{b\,z_1}\\
\\
\displaystyle
-
\left((\sqrt{\gamma_1}-\sqrt{\gamma_2})\,e^{2b\,z_1-(a-b)\,z_2}
+(\sqrt{\gamma_1}+\sqrt{\gamma_2})\,e^{2b\,z_1+(a-b)\,z_2}\right)
\,T_{12}\,\sqrt{\gamma_1}w'(0).
\end{array}
\tag {3.14}
$$

\noindent
Note that
$$\begin{array}{l}
\displaystyle
(\sqrt{\gamma_1}-\sqrt{\gamma_2})\,e^{2b\,z_1-(a-b)\,z_2}
+(\sqrt{\gamma_1}+\sqrt{\gamma_2})\,e^{2b\,z_1+(a-b)\,z_2}\\
\\
\displaystyle
=e^{2b\,z_1-(a-b)\,z_2}\,(\sqrt{\gamma_1}-\sqrt{\gamma_2})\left(1+O(e^{-2\,c\tau\,(a-b)/\sqrt{\gamma_2}})\right).
\end{array}
$$
Note also that Remark 1.1 gives $w'(0)=O(1)$ as $\tau\longrightarrow\infty$.
Using these, (3.11), (3.14), (2.4), (2.5) and the assumption $\gamma_1\not=\gamma_2$, we get the assertion of Lemma 3.1.

\noindent
$\Box$

Next we give a proof of Lemma 3.2.
Recalling (3.1) and (3.2), we get
$$\begin{array}{l}
\displaystyle
I_c(\tau)-I_c^0(\tau)\\
\\
\displaystyle
=-\gamma_1\,\Psi'(a)\epsilon(a)
+(\gamma_1-\gamma_2)\,\int_b^a\,\epsilon'(x)\,\Psi'(x)dx
-e^{-z^2\,T}\int_0^a\,u(\xi,T)\,\Psi(\xi)d\xi
\end{array}
\tag {3.15}
$$
where $\epsilon(y)=w(y)-y(x)$.  One knows
$$\displaystyle
\int_0^a u(\xi,T)\Psi(\xi)=O(1).
\tag {3.16}
$$

\noindent
The $\epsilon$ satisfies
$$\begin{array}{l}
\displaystyle
(\gamma\,\epsilon')'-z^2\,\epsilon=e^{-z^2\,T}u(x,T)\,\,\text{in}\,]0,\,a[,\\
\\
\displaystyle
\epsilon'(0)=\epsilon'(a)=0.
\end{array}
$$
Multiplying the both sides of the equation by $\overline\epsilon$ and integrating the resultant over
interval $]0,\,a[$, we have
$$\displaystyle
\int_0^a\gamma\vert\epsilon'\vert^2 dx+z^2\int_0^a\vert\epsilon\vert^2 dx
=-e^{-z^2\,T}\int_0^a\,u(\xi,T)\,\overline\epsilon(\xi)\,d\xi.
$$
Since
$$\displaystyle
z^2=\tau+i2\,c^2\,\tau^2\,\sqrt{1-\frac{1}{c^2\,\tau}},
$$
a standard argument yields,
as $\tau\longrightarrow\infty$
$$\displaystyle
\Vert\epsilon\Vert_{H^1(]0,\,a[)}=O(e^{-\tau\,T}).
$$

\noindent
A combination of this and the embedding $H^1(]0,\,a[)\subset C^0([0,\,a])$ gives the estimates
$$\begin{array}{l}
\displaystyle
-\gamma_1\,\Psi'(a)\epsilon(a)=
O(\tau\,\exp\,\left(-\tau\,(T+\frac{c\,a}{\sqrt{\gamma_1}})\right)),\\
\\
\displaystyle
\int_b^a\,\epsilon'(x)\,\Psi'(x)\,dx=
O(\tau\,\exp\,\left(-\tau\,(T+\frac{c\,b}{\sqrt{\gamma_1}})\right)).
\end{array}
\tag {3.17}
$$
Now we obtain from (3.15), (3.16) and (3.17) the assertion of Lemma 3.2.

\noindent
$\Box$

\section{Proof of Theorem 2.3.  Part 1.  Asymptotic behaviour of $w(a)$}

Define
$$\displaystyle
w(x)=w(x,\tau)=\int_0^T u(x,t)\,e^{-z^2\,t}\,dt,\,0<x<a.
$$
This $w(x)$ satisfies
$$\begin{array}{c}
\displaystyle
(\gamma\,w')'-z^2\,w=e^{-z^2\,T}\,u(x,T)\,\,\text{in}\,]0,\,a[,\\
\\
\displaystyle
\gamma_m\,w'(a)+\rho\,w(a)=0.
\end{array}
$$
In this section we study the asymptotic behaviour of $w(a)$ as $\tau\longrightarrow\infty$.
For the purpose it suffices to study the asymptotic behaviour of
the solution of the boundary value problem as $\tau\longrightarrow\infty$:
$$\begin{array}{c}
\displaystyle
(\gamma\,y')'-z^2y=0\,\,\text{in}\,]0,\,a[,\\
\\
\displaystyle
y'(0)=w'(0),\\
\\
\displaystyle
\gamma_m\,y'(a)+\rho\,y(a)=0.
\end{array}
$$

\noindent
This is because of
\proclaim{\noindent Lemma 4.1.}
Let $\rho\ge 0$.
The formula
$$\displaystyle
w(a)=y(a)+O(e^{-\tau\,T}),
$$
is valid.

\endproclaim

{\it\noindent Proof.}
Define
$$
\epsilon(x)=w(x)-y(x),\,\,0<x<a.
$$
This function satisfies
$$
\begin{array}{c}
\displaystyle
(\gamma\,\epsilon')'-z^2\,\epsilon=e^{-z^2\,T}u(x,T)\,\,\text{in}\,]0,\,a[,\\
\\
\displaystyle
\epsilon'(0)=0,\\
\\
\displaystyle
\gamma_m\,\epsilon'(a)+\rho\,\epsilon(a)=0.
\end{array}
$$
Hereafter a combination of a standard argument and the embedding
$H^1(]0,\,a[)\subset C^0([0,\,a])$ yields the desired estimate.

\noindent
$\Box$

For each $j=1,\cdots,m$ one can write
$$\displaystyle
y(x)=A_je^{x\cdot z_j}+B_je^{-xz_j},\,\,b_{j-1}<x<b_j.
$$
From the equation, it follows that
$$\begin{array}{c}
\displaystyle
y(b_j-0)=y_{j+1}(b_j+0),\\
\\
\displaystyle
\gamma_j\,y_j'(b_j-0)=\gamma_{j+1}\,y_{j+1}'(b_j+0).
\end{array}
$$
Therefore the coefficients $A_1,B_1,\cdots, A_m, B_m$ have to satisfy the system of equations:
$$\displaystyle
z_1(A_1-B_1)=w'(0);
\tag {4.1}
$$
for each $j=1,\cdots,m-1$
$$\begin{array}{c}
\displaystyle
A_je^{b_jz_j}+B_je^{-b_jz_j}
=A_{j+1}e^{b_jz_{j+1}}+B_{j+1}e^{-b_jz_{j+1}},\\
\\
\displaystyle
\gamma_j\,z_j(A_je^{b_jz_j}-B_je^{-b_jz_j})
=\gamma_{j+1}\,z_{j+1}(A_{j+1}e^{b_jz_{j+1}}-B_{j+1}e^{-b_jz_{j+1}});
\end{array}
\tag {4.2}
$$
$$\displaystyle
\gamma_m\,z_m(A_me^{az_m}-B_me^{-az_m})
+\rho(A_me^{az_m}+B_me^{-az_m})=0.
\tag {4.3}
$$

\noindent
Set
$$\begin{array}{lcr}
\displaystyle
X_j=\left(\begin{array}{c} A_j\,e^{b_jz_j}\\
\\
\displaystyle
B_j\,e^{-b_jz_j}\end{array}\right), &
\displaystyle
K_j=
\left(\begin{array}{cc} 1  & 1\\
\\
\displaystyle
\sqrt{\gamma_j} & -\sqrt{\gamma_j}
\end{array}
\right), &
\displaystyle
\alpha_j=\left(\begin{array}{cc}
1 & 0\\
\\
\displaystyle
0 & e^{-2(b_{j-1}-b_j)}\,z_j
\end{array}
\right).
\end{array}
$$

Then one can rewrite the equations (4.1), (4.2) and (4.3) in the
matrix form:
$$\displaystyle
\left(\begin{array}{cc} 1  & -e^{2b_1z_1}\end{array}\right)X_1=
\displaystyle\frac{w'(0)\,\sqrt{\gamma_1}}{z}e^{b_1\,z_1},
\tag {4.4}
$$
$$\displaystyle
e^{(b_j-b_{j+1})\,z_{j+1}}K_{j+1}\alpha_{j+1}X_{j+1}=K_jX_j,
\tag {4.5}
$$
$$\displaystyle
\left(\begin{array}{cc} \displaystyle(\sqrt{\gamma_m}+\frac{\rho}{z}) &
\displaystyle
-(\sqrt{\gamma_m}-\frac{\rho}{z})\end{array}\right)
X_m=0.
\tag {4.6}
$$

\noindent
Set
$$\displaystyle
L(z)=(K_1^{-1}\,K_2\,\alpha_2)\cdots(K_{m-1}^{-1}\,K_m\,\alpha_m).
$$
From (4.5) we have
$$\displaystyle
X_1=e^{\displaystyle\sum_{j=1}^{m-1}(b_j-b_{j+1})\,z_{j+1}}L(z)X_m.
$$

\noindent
Substituting this into (4.4), we obtain
$$\displaystyle
\left(\begin{array}{cc} 1  & -e^{2b_1z_1}\end{array}\right)L(z)X_m=
\displaystyle\frac{w'(0)\,\sqrt{\gamma_1}}{z}e^{b_1\,z_1}
e^{\displaystyle -\sum_{j=1}^{m-1}(b_j-b_{j+1})\,z_{j+1}}.
\tag {4.7}
$$

The equations (4.6) and (4.7) are equivalent to the equation
$$
\displaystyle
\left(\begin{array}{c}
\displaystyle\begin{array}{cc} \displaystyle 1 & -\frac{\displaystyle\sqrt{\gamma_m}-\frac{\rho}{z}}
{\displaystyle\sqrt{\gamma_m}+\frac{\rho}{z}}\end{array}\\
\\
\displaystyle
\left(\begin{array}{cc} 1 & -e^{2b_1\,z_1}\end{array}\right)L(z)
\end{array}
\right)\,X_m
=\displaystyle\frac{w'(0)\,\sqrt{\gamma_1}}{z}e^{b_1\,z_1}
e^{\displaystyle -\sum_{j=1}^{m-1}(b_j-b_{j+1})\,z_{j+1}}
\left(\begin{array}{c} 0\\
\\
\displaystyle 1\end{array}\right).
\tag {4.8}
$$
Solving equation (4.8), we obtain $X_m$.   The problem is the asymptotic behaviour
of $L(z)$ as $\tau\longrightarrow\infty$.

Define the {\it transmission coefficient} $T_{kl}$ and {\it reflection coefficient} $R_{kl}$
by the formula
$$\displaystyle
T_{kl}=\frac{2\sqrt{\gamma_k}}{\sqrt{\gamma_k}+\sqrt{\gamma_l}},
\,\,
R_{kl}
=\frac{\sqrt{\gamma_k}-\sqrt{\gamma_l}}
{\sqrt{\gamma_k}+\sqrt{\gamma_l}}.
$$
Set $\delta=\min_{j=1,\cdots,m}\,\{(b_j-b_{j-1})/\sqrt{\gamma_j}\}(>0)$.

\proclaim{\noindent Lemma 4.2.}
We have, as $\tau\longrightarrow\infty$
$$\displaystyle
L(z)=\frac{1}{\displaystyle T_{12}\cdots T_{m-1,\,m}}
\left(\begin{array}{cc}
1 & 0\\
\\
\displaystyle
R_{12} & 0
\end{array}
\right)
+O(e^{-2c\delta\tau}).
\tag {4.9}
$$

\endproclaim

{\it\noindent Proof.}
Using the expression
$$\displaystyle
\alpha_j=\left(\begin{array}{cc}
\displaystyle
1 & 0
\\
\\
\displaystyle
0 & 0
\end{array}
\right)
+O(e^{-2c\tau(b_j-b_{j-1})/\sqrt{\gamma_j}})
$$
and
$$\displaystyle
K_j^{-1}
=
\frac{1}{\displaystyle 2\sqrt{\gamma_j}}
\left(\begin{array}{cc}
\displaystyle
\sqrt{\gamma_j} & 1\\
\\
\displaystyle
\sqrt{\gamma_j} & -1
\end{array}\right),
$$
we have
$$\displaystyle
K_j\,\alpha_j\,K_j^{-1}
=
\frac{1}{\displaystyle 2\sqrt{\gamma_j}}
\left(\begin{array}{c}
\displaystyle
1\\
\\
\displaystyle
\sqrt{\gamma_j}
\end{array}
\right)
\left(\begin{array}{cc}
\displaystyle
\sqrt{\gamma_j} & 1
\end{array}
\right)
+O(e^{-2c\tau(b_j-b_{j-1})/\sqrt{\gamma_j}}).
$$
This gives
$$\begin{array}{l}
\displaystyle
(K_2\alpha_2K_2^{-1})\cdots(K_m\alpha_mK_m^{-1})\\
\\
\displaystyle
=(\frac{1}{2})^{m-1}
\frac{\displaystyle \Pi_{j=2}^{m-1}(\sqrt{\gamma_j}+\sqrt{\gamma_{j+1}})}
{\displaystyle\sqrt{\Pi_{j=2}^m\gamma_j}}
\left(\begin{array}{c} \displaystyle 1\\
\\
\displaystyle
\sqrt{\gamma_2}
\end{array}
\right)
\left(\begin{array}{cc}
\displaystyle
\sqrt{\gamma_m} & 1\end{array}\right)
+O(e^{-2c\delta\tau}).
\end{array}
$$

\noindent
Since
$$\displaystyle
K_1^{-1}\,
\left(\begin{array}{c} \displaystyle 1\\
\\
\displaystyle
\sqrt{\gamma_2}
\end{array}
\right)
\left(\begin{array}{cc}
\displaystyle
\sqrt{\gamma_m} & 1\end{array}\right)\,K_m
=\sqrt{\frac{\gamma_m}{\gamma_1}}\,(\sqrt{\gamma_1}+\sqrt{\gamma_2})\,
\left(\begin{array}{cc}
\displaystyle
1 & 0
\\
\\
\displaystyle
\frac{\sqrt{\gamma_1}-\sqrt{\gamma_2}}
{\sqrt{\gamma_1}+\sqrt{\gamma_2}} & 0
\end{array}
\right),
$$
we obtain
$$\displaystyle
L(z)
=(\frac{1}{2})^{m-1}
\frac{\displaystyle \Pi_{j=1}^{m-1}(\sqrt{\gamma_j}+\sqrt{\gamma_{j+1}})}
{\displaystyle\sqrt{\Pi_{j=1}^{m-1}\gamma_j}}
\left(\begin{array}{cc}
\displaystyle
1 & 0
\\
\\
\displaystyle
\frac{\sqrt{\gamma_1}-\sqrt{\gamma_2}}{\sqrt{\gamma_1}+\sqrt{\gamma_2}} & 0
\end{array}
\right)
+O(e^{-2c\delta\tau}).
$$
This is nothing but (4.9).

\noindent
$\Box$

\noindent
As a direct consequence of (4.9) we have
$$\begin{array}{l}
\displaystyle
\left(\begin{array}{c}
\displaystyle\begin{array}{cc} \displaystyle 1 & -\frac{\displaystyle\sqrt{\gamma_m}-\frac{\rho}{z}}
{\displaystyle\sqrt{\gamma_m}+\frac{\rho}{z}}\end{array}\\
\\
\displaystyle
\left(\begin{array}{cc} 1 & -e^{2b_1\,z_1}\end{array}\right)L(z)
\end{array}
\right)^{-1}\\
\\
\displaystyle
=
\left(\begin{array}{lr}
\displaystyle 0 & T_{12}\cdots T_{m-1,m}\\
\\
\displaystyle -1 & T_{12}\cdots T_{m-1,m}\\
\\
\end{array}
\right)
+\frac{2\rho}{z\,\sqrt{\gamma_m}}\,
\left(\begin{array}{lr}
\displaystyle 0 & 0\\
\\
\displaystyle
-1 & T_{12}\cdots T_{m-1,m}\end{array}
\right)
+O\left(\frac{1}{\tau^2}\right)
\end{array}
$$
and (4.8) therefore yields
$$
\displaystyle
X_m=
\displaystyle\frac{w'(0)\,\sqrt{\gamma_1}}{z}e^{b_1\,z_1}
e^{\displaystyle -\sum_{j=1}^{m-1}(b_j-b_{j+1})\,z_{j+1}}
T_{12}\cdots T_{m-1,\,m}
\{
\left(\begin{array}{c} \displaystyle 1\\
\\
\displaystyle 1+\frac{2\rho}{z\,\sqrt{\gamma_m}}\end{array}\right)
+O\left(\frac{1}{\tau^2}\right)\}
$$
as $\tau\longrightarrow\infty$.

Define
$$\displaystyle
\varphi_j=\sum_{l=1}^{j-1}b_l\,(z_l-z_{l+1}), \,\,j=2,\cdots,m
$$
and
$$
\varphi(x)=\left\{
\begin{array}{lr}
\displaystyle x\,z_1,&\quad \text{if $-\infty<x<b_1$,}\\
\\
\displaystyle x\,z_2+\varphi_2,&\quad \text{if $b_1<x<b_2$,}\\
\\
\vdots\\
\\
\displaystyle x\,z_m+\varphi_m, &\quad\text{if $b_{m-1}<x<\infty$.}
\end{array}
\right.
$$
The function $\varphi$ has the unique continuous extension to the
whole real line since $\varphi(b_j-0)=\varphi(b_j+0)$ for each
$j=1,\cdots,m-1$.  We denote the extension by $\varphi$ again.

\noindent
Since $y(a)=(X_m)_1+(X_m)_2$, $a=b_m$ and
$$\begin{array}{l}
\displaystyle
b_1\,\,z_1-\sum_{j=1}^{m-1}(b_j-b_{j+1})\,z_{j+1}\\
\\
\displaystyle
=b_1\,z_1-(b_1-b_2)\,z_2-(b_2-b_3)\,z_3-\cdots
-(b_{m-1}-b_m)\,z_m\\
\\
\displaystyle
=b_1(z_1-z_2)+b_2\,(z_2-z_3)+\cdots+b_{m-1}\,(z_{m-1}-z_m)+b_m\,z_m
=\varphi(b_m),
\end{array}
$$
we obtain the formula
$$
\displaystyle
y(a)
=\frac{2w'(0)\,\sqrt{\gamma_1}}{z}\,e^{\varphi(a)}\,
T_{12}\cdots T_{m-1,\,m}
\{1+\frac{\rho}{z\,\sqrt{\gamma_m}}+O\left(\frac{1}{\tau^2}\right)\}.
\tag {4.10}
$$
Note that
$$\varphi(a)=-z\,\left(\frac{b_1}{\sqrt{\gamma_1}}+\sum_{j=1}^{m-1}\,\frac{b_{j+1}-b_j}
{\sqrt{\gamma_{j+1}}}\right).
\tag {4.11}
$$

\section{Proof of Theorem 2.3.  Part 2.  Asymptotic behaviour of $\Psi(a)$}

Integration by parts yields that
the function $\Psi$ is the (weak) solution of the equation
$(\tilde{\gamma}\,y')'-z^2y=0$ in $\Bbb R$ if and only if
$\Psi$ satisfies, for each $j=1,\cdots,m-1$
$$\begin{array}{c}
\displaystyle
\Psi(b_j-0)=\Psi(b_j+0),\\
\\
\displaystyle
\gamma_j\,\Psi'(b_j-0)=\gamma_{j+1}\,\Psi'(b_j+0).
\end{array}
$$

This yields the system of equations for $B_1, A_2,B_2,\cdots, B_{m-1},
A_m$:
$$\displaystyle
e^{(b_j-b_{j+1})\,z_{j+1}}K_{j+1}\alpha_{j+1}Y_{j+1}=K_jY_j,\,\,
j=1,\cdots,m-1
\tag {5.1}
$$
where
$$\begin{array}{lcr}
\displaystyle
Y_1=\left(\begin{array}{c}\displaystyle e^{b_1\,z_1}\\
\\
\displaystyle
B_1 e^{-b_1\,z_1}
\end{array}
\right),  &
Y_j=\left(\begin{array}{c}
\displaystyle
A_j e^{b_j\,z_j}\\
\\
\displaystyle
B_j e^{-b_j\,z_j}
\end{array}
\right),  &
Y_m=A_m\,e^{b_m\,z_m}\left(\begin{array}{c}\displaystyle
1\\
\\
\displaystyle
0
\end{array}
\right).
\end{array}
$$

From (5.2) one has
$$\begin{array}{l}
\displaystyle
Y_1=e^{\displaystyle\sum_{j=1}^{m-1}(b_j-b_{j+1})\,z_{j+1}}
L(z)Y_m\\
\\
\displaystyle
=A_m
\,e^{b_m\,z_m}\,e^{\displaystyle\sum_{j=1}^{m-1}(b_j-b_{j+1})\,z_{j+1}}L(z)\left(\begin{array}{c}
\displaystyle 1\\
\\
\displaystyle
0
\end{array}
\right).
\end{array}
\tag {5.3}
$$

\noindent
Set $\delta=\min_{j=1,\cdots,m}\,\{(b_j-b_{j-1})/\sqrt{\gamma_j}\}(>0)$.

\noindent
The following indicates the meaning of the coefficients $T_{kl}$
and $R_{kl}$.

\proclaim{\noindent Lemma 5.1.}
For sufficiently large $\tau$ the system of equations (5.1) is uniquely solvable
and, as $\tau\longrightarrow\infty$ the formulae
$$\displaystyle
A_j\,e^{-b_1\,z_1}\,e^{b_j\,z_j}\,
e^{\displaystyle
\sum_{l=1}^{j-1}(b_{l}-b_{l+1})\,z_{l+1}}
=T_{12}\cdots T_{j-1,\,j}+O(e^{-2c\delta\tau})
\tag {5.4}
$$
and
$$
\displaystyle
B_j\,e^{-b_1\,z_1}\,e^{-b_j\,z_j}\,
e^{\displaystyle
\sum_{l=1}^{j-1}(b_{l}-b_{l+1})\,z_{l+1}}
=R_{j,\,j+1}\,T_{12}\cdots T_{j-1,\,j}
+O(e^{-2c\delta\tau}),
\tag {5.5}
$$
are valid.

\endproclaim

{\it\noindent Proof.}
Using (4.9) and (5.3), we find $A_m$ for a sufficient large $\tau$
and obtain (5.4) for $j=m$.
Next (5.1) for $j=m-1$ yields that $A_{m-1}$ and $B_{m-1}$ are uniquely determined
and that the formulae (5.4) and (5.5) for $j=m-1$.
For general $j$ we make use of the recurrence formulae
$$\begin{array}{l}
\displaystyle
A_j=\frac{1}{T_{j,\,j+1}}A_{j+1}\,e^{b_j\,(z_{j+1}-z_j)}
+\frac{R_{j,\,j+1}}{T_{j,\,j+1}}B_{j+1}\,e^{-b_j\,(z_{j+1}+z_j)},\\
\\
\displaystyle
B_j=\frac{R_{j,\,j+1}}{T_{j,\,j+1}}A_{j+1}\,e^{b_j(z_{j+1}+z_j)}
+\frac{1}{T_{j,\,j+1}}B_{j+1}\,e^{-b_j\,(z_{j+1}-z_j)},
\end{array}
$$
which are equivalent to (5.1).
Note that we set $A_1=1$ and $B_m=0$.

\noindent
$\Box$

\noindent
From Lemma 5.1 one has, as $\tau\longrightarrow\infty$
$$\begin{array}{c}
\displaystyle
A_2\,e^{-\varphi_2}=T_{12}+O(e^{-2c\delta\tau}),\\
\\
\displaystyle
A_3\,e^{-\varphi_3}=T_{12}\,T_{23}+O(e^{-2c\delta\tau}),\\
\\
\displaystyle
\vdots
\\
\\
\displaystyle
A_m\,e^{-\varphi_m}=T_{12}\,T_{23}\cdots T_{m-1,m}+O(e^{-2c\delta\tau});
\\
\\
\displaystyle
B_1\,e^{-2b_1\,z_1}=R_{12}+O(e^{-2c\delta\tau}),\\
\\
\displaystyle
B_2\,e^{-2b_2\,z_2}\,e^{-\varphi_2}
=T_{12}\,R_{23}+O(e^{-2c\delta\tau}),\\
\\
\displaystyle
\vdots\\
\\
\displaystyle
B_{m-1}\,e^{-2b_{m-1}\,z_{m-1}}
\,e^{-\varphi_{m-1}}
=T_{12}\,T_{23}\cdots T_{m-2,\,m-1}R_{m-1,\,m}+O(e^{-2c\delta\tau}).
\end{array}
$$

\noindent
Moreover, it follows that
$$
\frac{B_j\, e^{-x\,z_j}}
{A_j\,e^{x\,z_j}}=O(e^{-2c\tau(b_j-x)/\sqrt{\gamma_j}}),\,\,x<b_j.
$$
This gives, as $\tau\longrightarrow\infty$ the asymptotic formula of $\Psi$:
$$
\Psi(x)\sim\left\{
\begin{array}{lr}
\displaystyle e^{\varphi(x)},&\quad \text{if $-\infty<x<b_1$,}\\
\\
\displaystyle T_{12}\,e^{\varphi(x)},&\quad \text{if $b_1<x<b_2$,}\\
\\
\vdots\\
\\
\displaystyle T_{12}\cdots T_{m-1,\,m}\,e^{\varphi(x)}, &\quad\text{if $b_{m-1}<x<\infty$.}
\end{array}
\right.
$$

\noindent
The formula (5.4) for $j=m$ gives
$$\displaystyle
\Psi(a)\,e^{-\varphi(a)}
=T_{12}\,\cdots T_{m-1,\,m}
+O(e^{-2c\delta\tau}).
\tag {5.6}
$$

\noindent The following estimates are a direct corollary of (5.4)
and (5.5):
$$\begin{array}{l}
\displaystyle
\vert A_j\vert
=O
\left(e^{-c\,b_1\tau/\sqrt{\gamma_1}}\,e^{c\,b_j\tau/\sqrt{\gamma_j}}\,
e^{\displaystyle -c\sum_{l=1}^{j-1}(b_{l+1}-b_l)\,\tau/\sqrt{\gamma_{l+1}}}\right);\\
\\
\displaystyle
\vert B_j\vert
=O\left(e^{-c\,b_1\tau/\sqrt{\gamma_1}}\,e^{-c\,b_j\tau/\sqrt{\gamma_j}}\,
e^{\displaystyle -c\,\sum_{l=1}^{j-1}(b_{l+1}-b_l)\,\tau/\sqrt{\gamma_{l+1}}}\right).
\end{array}
$$

\noindent
Applying these estimates to the expression
$$\displaystyle
\int_{b_{j-1}}^{b_j}u(\xi,T)\Psi(\xi)d\xi
=A_j\int_{b_{j-1}}^{b_j}u(\xi,T)e^{\xi\,z_j}d\xi
+B_j\int_{b_{j-1}}^{b_j}u(\xi,T)e^{-\xi\,z_j}d\xi,
$$
we obtain the estimate
$$\displaystyle
\int_0^au(\xi,T)\Psi(\xi)d\xi
=O(1).
\tag {5.7}
$$

Recalling the definition of $w$, the equation $\Psi'a)=z_m\Psi(a)$ and using integration by parts we
have
$$\begin{array}{l}
\displaystyle
I_{c}(\tau)
=\gamma_1\,w'(0)\Psi(0)
-w(0)\,\gamma_1\Psi'(0)\\
\\
\displaystyle
=-(\rho\,+\gamma_m\,z_m)\Psi(a)w(a)
-e^{-z^2T}\int_0^a u(\xi,T)\Psi(\xi)d\xi.
\end{array}
$$

\noindent
Then Lemma 4.1, (4.10), (5.6) and (5.7) yield
the asymptotic formula for the indicator function.

\proclaim{\noindent Proposition 5.2.}
As $\tau\longrightarrow\infty$ the formula
$$\begin{array}{l}
\displaystyle
I_{c}(\tau)\,e^{-2\,\varphi(a)}
=-2\,\sqrt{\gamma_1\,\gamma_m}\,w'(0)\,(T_{12}\cdots T_{m-1,\,m})^2
\{1+\frac{2\rho}{z\,\sqrt{\gamma_m}}+O(\frac{1}{\tau^2})\}
\\
\\
\displaystyle
+O(\tau e^{\displaystyle -\tau\{T-
2c\,(
\frac{b_1}{\sqrt{\gamma_1}}+
\sum_{j=1}^{m-1}\frac{b_{j+1}-b_j}{\sqrt{\gamma_{j+1}}}
)\}}),
\end{array}
$$
is valid.

\endproclaim

\noindent

Theorem 2.3 is an immediate corollary of (4.11), Proposition 5.2,
the expression
$$\displaystyle
w'(0)=\int_0^T e^{-z^2\,t}u_x(0,t)dt,
$$
the assumption (2.7) and the fact $w'(0)=O(1)$ as $\tau\longrightarrow\infty$.

{\bf\noindent Remark 5.1.}
Once $\varphi(a)$ is known (see (4.11)), then one immediately obtains the extraction formula of $\rho$:
$$
\displaystyle
\lim_{\tau\longrightarrow\infty}\left(\frac{I_c(\tau)\,e^{-2\,\varphi(a)}}{2\sqrt{\gamma_1\,\gamma_m}\,w'(0)\,
(T_{12}\cdots T_{m-1,\,m})^2}+1\right)\,z\,\sqrt{\gamma_m}=-2\rho.
$$

\section{Proof of Theorem 2.5}

Define
$$\displaystyle
w(x,\tau)=\int_0^T e^{-z^2\,t}\,u(x,t)dt,\,\,0<x<a.
$$
This $w$ satisfies
$$\begin{array}{c}
\displaystyle
(\gamma\,w')'-z^2\,w=e^{-z^2\,T}u(x,T),\,\,\text{in}\,]0,\,a[,\\
\\
\displaystyle
\gamma(a)\,w'(a)+\rho\,w(a)=0.
\end{array}
$$
Integration by parts gives the expression
$$
I(z)
=-w(a,z)(\gamma(a)\,\Psi'(a;z,M)+\rho\,\Psi(a;z,M))
-\int_0^a\,e^{-z^2\,T}u(x,T)\Psi(x;z,M)dx.
\tag {6.1}
$$
Proposition 2.4 gives
$$\displaystyle
\gamma(a)\,\Psi'(a;z,M)+\rho\,\Psi(a;z,M)
=\frac{z}{\gamma(a)^{3/4}}e^{K_a\,z\,\pi}
\{1+O\left(\frac{1}{\vert z\vert}\right)\}.
\tag {6.2}
$$
Therefore it suffices to study the asymptotic behaviour of
$w(a,z)$ as $\vert z\vert\longrightarrow\infty$. For the purpose
first we study the asymptotic behaviour of the unique solution of
the boundary value problem as $\vert z\vert\longrightarrow\infty$:
$$\begin{array}{c}
\displaystyle
(\gamma y')'-z^2y=0\,\,\text{in}\,]0,\,a[,\\
\\
\displaystyle
y'(0)=1,\\
\\
\displaystyle
\gamma(a)y'(a)+\rho\,y(a)=0.
\end{array}
$$
We make use of the Liouville transform.
Define
$$\displaystyle
K_a=\frac{1}{\pi}\int_0^a\frac{dx}{\sqrt{\gamma(x)}}
$$
and
$$
\displaystyle
s(x;a)=\frac{1}{K_a}\int_0^x\frac{dx}{\sqrt{\gamma(x)}},\,\,0\le x\le a.
$$
Denote by $x(s;a),\,0\le s\le\pi$ the inverse of the function $s=s(x;a)$.
Set
$$\displaystyle
\tilde{y}(s)=K_a\gamma(0)^{1/4}\gamma(x(s;a))^{1/4}y(x(s;a)).
\tag {6.3}
$$
Then this $\tilde{y}$ satisfies
$$\begin{array}{c}
\displaystyle
\tilde{y}''-(K_a^2z^2+g_a(s))\tilde{y}=0\,\,\text{in}\,]0,\,\pi[,\\
\\
\displaystyle
\tilde{y}'(0)-h_a\,\tilde{y}(0)=1,\\
\\
\displaystyle
\tilde{y}'(\pi)+H_a\,\tilde{y}(\pi)=0
\end{array}
$$
where
$$\displaystyle
h_a=\frac{\gamma'(0)}{4K_a\,\gamma(0)^{3/2}},\,\,
H_a=\frac{4\rho-\gamma'(a)}
{4K_a\,\gamma(a)^{3/2}}.
$$
and
$$\displaystyle
f_a(s)=\gamma(x(s;a))^{1/4},\,\,
g_a(s)=\frac{f_a''(s)}{f_a(s)}.
$$

\proclaim{\noindent Lemma 6.1.}
As $\vert z\vert\longrightarrow\infty$ we have
$$\displaystyle
\tilde{y}(\pi)=\frac{2}{K_a\,z}e^{K_a\,z\,\pi}\{1+O\left(\frac{1}{\vert z\vert}\right)\}.
\tag {6.4}
$$
\endproclaim

{\it\noindent Proof.}
It is easy to see that $\tilde{y}$ has the expression
$$\displaystyle
\tilde{y}(s)=\frac{\varphi(s)}{\varphi'(0)-h_a\varphi(0)}
$$
where $\varphi$ is the unique solution of the initial value
problem:
$$\begin{array}{c}
\displaystyle
\varphi''-(K_a^2\,z^2+g(s))\varphi=0\,\,\text{in}\,]0,\,\pi[,\\
\\
\displaystyle
\varphi(\pi)=1,\\
\\
\displaystyle
\varphi'(\pi)=-H_a.
\end{array}
$$
Here we cite a known important fact:
there exists a fundamental system of solutions $e_1(x,z)$ and $e_2(x,z)$
of equation $y''-(K_a^2\,z^2+g_a(s))y=0$ in $]0,\,\pi[$ such that
as $\vert z\vert\longrightarrow\infty$, $\text{Re}\,z\le 0$, uniformly in $x$,
$$\begin{array}{c}
\displaystyle
e_1(s,z)=e^{K_a\,z\,s}
\{1+O\left(\frac{1}{\vert z\vert}\right)\},\\
\\
\displaystyle
e_1'(s,z)=K_a\,z\,e^{K_a\,z\,s}\{1+O\left(\frac{1}{\vert z\vert}\right)\},\\
\\
\displaystyle
e_2(s,z)=e^{-K_a\,z\,s}\{1+O\left(\frac{1}{\vert z\vert}\right)\},\\
\\
\displaystyle
e_2'(s,z)=-K_a\,z\,e^{-K_a\,z\,s}\{1+O\left(\frac{1}{\vert z\vert}\right)\}
\end{array}
$$
See again Theorem 1 of p. 48 in \cite{N}.

\noindent
Thus one can write
$$
\varphi(s)=A_1\,e_1(s,z)+A_2\,e_2(s,z)
$$
where $A_1$ and $A_2$ are constants. The initial conditions on
$\varphi$ yields
$$\begin{array}{c}
\displaystyle
A_1=\frac{1}{W(e_1(\,\cdot\,,z),e_2(\,\cdot\,,z))(\pi)}(e_2'(\pi,z)+e_2(\pi,z)H_a),\\
\\
\displaystyle
A_2=-\frac{1}{
W(e_1(\,\cdot\,,z),e_2(\,\cdot\,,z))(\pi)
}(e_1'(\pi,z)+e_1(\pi,z)H_a).
\end{array}
$$
The asymptotic behaviour of $e_1$ and $e_2$ yields
$$\displaystyle
W(e_1(\,\cdot\,,z),e_2(\,\cdot\,,z))(\pi)
=-2K_a\,z
\{1+O\left(\frac{1}{\vert z\vert}\right)\}
$$
and also
$$\displaystyle
A_1=\frac{1}{2}e^{-K_a\,z\,\pi}
\{1+O\left(\frac{1}{\vert z\vert}\right)\},\,\,
A_2=\frac{1}{2}e^{K_a\,z\,\pi}
\{1+O\left(\frac{1}{\vert z\vert}\right)\}.
$$
From these we obtain
$$\begin{array}{c}
\displaystyle
\varphi'(0)-h_a\,\varphi(0)
=\frac{K_a\,z}{2}e^{-K_a\,z\pi}
\{1+O\left(\frac{1}{\vert z\vert}\right)\},\\
\\
\displaystyle
\varphi(\pi)=1+O\left(\frac{1}{\vert z\vert}\right).
\end{array}
$$
This yields (6.4).

\noindent
$\Box$

\noindent
A combination of (6.3) and (6.4) gives
$$\displaystyle
y(a)=\frac{2\,e^{K_a\,z\,\pi}}{K_a^2\,\gamma(0)^{1/4}\,\gamma(a)^{1/4}z}
\{1+O\left(\frac{1}{\vert z\vert}\right)\}.
\tag {6.5}
$$

Using a standard argument, we have
$$\displaystyle
\frac{w(x,z)}
{w'(0,z)}=y(x)+O(e^{\displaystyle -\text{Re}\,z^2\,T}).
\tag {6.6}
$$
A combination of (6.5) and (6.6) yields
$$\displaystyle
w(a,z)=\frac{2w'(0,z)}{K_a^2\,\gamma(0)^{1/4}\,\gamma(a)^{1/4}}
\frac{e^{K_a\,z\,\pi}}{z}
\{1+O\left(\frac{1}{\vert z\vert}\right)\}
+O(e^{\displaystyle -\text{Re}\,z^2\,T}w'(0,z)).
\tag {6.7}
$$
Now from (6.1), (6.2) and (6.7) we obtain the crucial formula:
$$\begin{array}{c}
\displaystyle
I(z)
=-\frac{2\,w'(0,z)}
{K_a^2\,\gamma(0)^{1/4}}
e^{2\,K_a\,z\,\pi}\{1+O\left(\frac{1}{\vert z\vert}\right)\}
\\
\\
\displaystyle
+O\left(e^{\displaystyle -\text{Re}\,z^2\,T}\vert w'(0,z)\vert\,\vert z\vert\,
e^{\displaystyle K_a\,\text{Re}\,z\pi}\right)
+O\left(e^{\displaystyle -\text{Re}\,z^2\,T}\right).
\end{array}
\tag {6.8}
$$
Then Theorem 2.5 follows from this formula (6.8),
for the choices $z=-c\tau(1+i\sqrt{1-1/c^2\tau})$ (case (a))
and $z=-\tau$ (case (b)) and the expression
$$
\displaystyle w'(0,z)=\int_0^T e^{-z^2\,t}u_x(0,t)dt.
$$
More precisely, both (2.11) in case (a) and (2.13) in case (b) ensure
$$\displaystyle
C\tau^{\mu}\le\vert w'(0,z)\vert,\,\,\tau\ge\tau^0.
$$
From Remark 1.1 one has $w'(0,z)=O(1)$ as $\tau\longrightarrow\infty$.
These yields the estimate: for suitable positive constants $C_1$, $C_2$ and for all $\tau>>1$
$$\displaystyle
C_1\tau^{\mu}\le\vert I(z)\vert e^{\displaystyle -2\,K_a\,\text{Re}\,z\pi}\le C_2.
$$
Note that (2.10) is essential in the case (a).  This completes the proof of Theorem 2.5.

$$\quad$$

\centerline{{\bf Acknowledgement}}

This research was partially supported by Grant-in-Aid for
Scientific Research (C)(No.  18540160) of Japan  Society for
the Promotion of Science.

\vskip1cm
\noindent
e-mail address

ikehata@math.sci.gunma-u.ac.jp
\end{document}